\documentclass[11pt]{amsart}

\title{Compactification for Essentially Finite-Type Maps}
\date{}

   \usepackage{amsthm}
   \usepackage{amsmath}
   \usepackage{amssymb}
   \usepackage{amsfonts}
   \usepackage{amscd}
   \usepackage[mathscr]{eucal}
   \usepackage{verbatim}
   \usepackage{graphics}

\input xy
\xyoption{all}

\newtheorem{athm}{Theorem}[section]
\newtheorem{alem}[athm]{Lemma}
\newtheorem{aprop}[athm]{Proposition}

\theoremstyle{definition}
 \newtheorem{arem}[athm]{Remark}
 \newtheorem{arems}[athm]{Remarks}
 \newtheorem{adefi}[athm]{Definition}

\numberwithin{equation}{athm}

\newcommand{\sss}{\stepcounter{athm}\subsection{}}

\newcommand{\ssss}{\subsubsection{}}
\newcommand{\EQAL}[1]%
{\,\begin{picture}(#1,0)%
\put(0,3){\line(1,0){#1}}%
\put(0,1){\line(1,0){#1}}%
\end{picture}\,}%

\newcommand{\vlto}[1]%
{\,\begin{picture}(#1,3)%
\put(0,2){\vector(1,0){#1}}%
\end{picture}\,}%

\newcommand{\vllarrow}[1]%
{\,\begin{picture}(#1,3)%
\put(#1,2){\vector(-1,0){#1}}%
\end{picture}\,}%

\newcommand{\dirlm}[1]%
{ {\lim\hskip-1.58em\lower.65ex \hbox{$ {}_{\stackrel{\lower1ex\hbox
   {$\scriptstyle -\!\!\<\longrightarrow$}}{ ^{#1} } } $}} \:}

\newcommand{\subdirlm}[1]%
  { {\lim\hskip-1.5em\lower.6ex \hbox{$ {}_{\stackrel{\lower1ex\hbox
  {$\scriptstyle\longrightarrow$}}{ ^{#1} } } $} } \:}

\newcommand{\inlm}[1]%
  { {\lim\hskip-1.58em\lower.65ex  \hbox{$ {}_{\stackrel{\lower1ex\hbox
  {$\scriptstyle \longleftarrow\!\!\<-$}  }{ ^{#1} }  } $}} \:}

\def\>{\mspace {1mu}}
\def\<{\mspace{-1mu}}
\def\({{\textup(}}
\def\){{\textup)}}

\newcommand{\X}{{\mathscr X}}

\newcommand{\Y}{{\mathscr Y}}
\newcommand{\Z}{{\mathscr Z}}
\newcommand{\V}{{\mathscr V}}

\newcommand{\W}{{\mathscr W}}
\newcommand{\I}{{\mathscr I}}
\newcommand{\J}{{\mathscr J}}

\newcommand{\sS}{{\mathscr S}}
\newcommand{\sT}{{\mathscr T}}
\newcommand{\sU}{{\mathscr U}}
\newcommand{\sV}{{\mathscr V}}

\newcommand{\eC}{{\mathscr C}}

\newcommand{\sfC}{{\mathsf C}}
\newcommand{\sfI}{{\mathsf I}}
\newcommand{\sfF}{{\mathsf F}}
\newcommand{\sfO}{{\mathsf O}}
\newcommand{\sfP}{{\mathsf P}}
\newcommand{\sfQ}{{\mathsf Q}}
\newcommand{\sfR}{{\mathsf R}}

\newcommand{\Cfr}{{\mathfrak C}}

\newcommand{\A}{{\mathcal A}}
\newcommand{\B}{{\mathcal B}}
\newcommand{\C}{{\mathcal C}}
\newcommand{\cD}{{\mathsf D}}
\newcommand{\E}{{\mathcal E}}
\newcommand{\F}{{\mathcal F}}

\renewcommand{\H}{{\mathcal H}}

\newcommand{\cI}{{\mathcal I}}

\newcommand{\cS}{{\mathcal S}}

\newcommand{\cFb}{{\mathcal F}^{\bullet}}

\newcommand{\cMb}{{\mathcal M}^{\bullet}}

\newcommand{\cO}{{\mathcal O}}
\newcommand{\OX}{\cO_{X}}
\newcommand{\OY}{\cO_{Y}}
\newcommand{\OU}{\cO_{U}}

\newcommand{\OZ}{\cO_{Z}}
\newcommand{\OXx}{\cO_{X,x}}
\newcommand{\OYy}{\cO_{Y,y}}

\newcommand{\OZz}{\cO_{Z,z}}

\newcommand{\Spec}{{\mathrm {Spec}}}

\newcommand{\Bl}{{\mathrm {Bl}}}

\newcommand{\cont}{{\mathrm {c}}}
\newcommand{\rh}{{\mathrm {h}}}
\newcommand{\rv}{{\mathrm {v}}}

\newcommand{\bC}{{\mathbf C}}
\newcommand{\D}{{\mathbf D}}

\newcommand{\Dqc}{\D_{\mkern-1.5mu\mathrm {qc}}}

\newcommand{\Dqct}{\D_{\mkern-1.5mu\mathrm{qct}}}

\newcommand{\Dc}{\D_{\mkern-1.5mu\mathrm c}}

\newcommand{\bt}{{\bf t}}
\newcommand{\bs}{{\bf s}}

\newcommand{\boldC}{\boldsymbol{C}}
\newcommand{\bbeta}{\boldsymbol{\beta}}
\newcommand{\bphi}{\boldsymbol{\phi}}

\newcommand{\cfr}{\mathfrak c}

\newcommand{\lfr}{\mathfrak l}

\newcommand{\pfr}{\mathfrak p}
\newcommand{\qfr}{\mathfrak q}
\newcommand{\rfr}{\mathfrak r}
\newcommand{\sfr}{\mathfrak s}

\newcommand{\Se}{\textbf{\textup{S}}_{\textbf{\textup{e}}}}
\newcommand{\Sf}{\textbf{\textup{S}}_{\textbf{\textup{f}}}}

\newcommand{\R}{{\mathbf R}}

\newcommand{\bL}{{\mathbf L}}

\newcommand{\Hom}{{\mathrm {Hom}}}
\newcommand{\Homc}{{\mathrm {Hom}}^{\cont}}

\newcommand{\sHom}{\H{om}}

\newcommand{\iGp}[1]{{\varGamma_{\<\!#1}'}}
\newcommand{\iG}[1]{{\varGamma_{\<\!#1}^{\phantom\prime}}}

\newcommand{\set}{\!:=}

\newcommand{\Lra}{\Longrightarrow}

\newcommand{\onto}{\twoheadrightarrow}

\newcommand{\xto}[1]{\xrightarrow{#1}}
\newcommand{\xgets}[1]{\xleftarrow{#1}}

\newcommand{\ssbox}{{\scriptscriptstyle{\Box}}}
\newcommand{\bxt}{\scriptscriptstyle{\boxtimes}}
\newcommand{\sbsq}{{\scriptstyle{\blacksquare}}}

\newcommand{\wt}[1]{{\widetilde{#1}}}
\newcommand{\ov}[1]{{\overline{#1}}}

\newcommand{\wtcD}{\wt{\cD}}

\newcommand{\wtDqc}{\wt{\D}_{\mkern-1.5mu\mathrm {qc}}}

\newcommand{\sbeta}{\textup{\ss}}

\newcommand{\dg}[1]{\ensuremath{\dagger_{\<\raisebox{-.2ex}{\scriptsize{#1}}}}}
\newcommand{\oneD}[1]{\smash{\mathbf{1}_{\<\<%
\raisebox{-.1ex}{$\scriptstyle\cD_{\<\<%
\raisebox{-.1ex}{$\scriptscriptstyle{#1}$}}$}}}}
\newcommand{\xyoneD}[1]{\;\;\mathbf{1}_{\<\<%
\raisebox{-.1ex}{$\scriptstyle\cD_{\<\<%
\raisebox{-.1ex}{$\scriptscriptstyle{#1}$}}$}}}

\newcommand{\iso}%
{{\mkern8mu\longrightarrow \mkern-25.5mu{}^\sim\mkern17mu}}
\newcommand{\osi}%
{{\mkern8mu\longleftarrow \mkern-24.5mu{}^\sim\mkern16mu}}
\def\Otimes{\underset
  {\vbox to 0pt {\vskip-1ex\hbox{$\scriptscriptstyle=$}\vss}}
    {\otimes}\vadjust{\kern.4pt}}

\newcommand{\smcirc}%
  {{\raise.15ex\hbox to.7em{$\hss \scriptstyle\circ\hss$}}} 

\begin{document}

\author{Suresh Nayak}
\address{Chennai Mathematical Institute \\ Plot H1, SIPCOT IT Park,
Siruseri-603103, INDIA}
\address{\textit{Current address:} 
Avery Hall, University of Nebraska, 
Lincoln, NE-68588, USA}
\email {snayak@cmi.ac.in, snayak@math.purdue.edu}
\thanks{This research has been carried out during the author's
visit to Purdue University and the University of Nebraska at Lincoln.
The author is grateful to them for their hospitality.}

\begin{abstract}
We show that any
separated essentially finite-type 
map~$f$ of noetherian schemes globally
factors as $f = hi$ where~$i$ is an injective 
localization map and~$h$ a separated 
finite-type map. In particular, via Nagata's compactification
theorem, $h$ can be chosen to be proper. 
We apply these results to Grothendieck duality. We also
obtain other factorization results and 
provide essentialized versions of many general
results such as Zariski's Main Theorem, 
Chow's Lemma, and blow-up descriptions 
of birational maps.
\end{abstract}

\maketitle

\renewcommand{\theenumi}{\roman{enumi}}

\section{Introduction}
\label{sec:Intro}
Let us denote by $\Sf$ the category of separated finite-type 
morphisms of noetherian
schemes. Nagata's compactification theorem 
(\cite{Ng}, \cite{Con2}, \cite{Lu}) states 
that any $\Sf$-map $f \colon X \to S$ factors as 
\smash{$X \xto{\;i\;} \ov{X} \xto{\;g\;} S$} where~$i$ is an 
open immersion and~$g$ is a proper map. 
One of its applications lies in Grothendieck 
duality where the only known approach to defining the 
twisted-inverse-image pseudofunctor~$(-)^!$ over all of~$\Sf$ 
relies on the compactification theorem. 
This approach, pioneered by Deligne and Verdier
(\cite{De}, \cite{V}),
and developed further
by Lipman, Neeman, Sastry and others
(\cite{Lip-LNM}, \cite{DFS}, \cite{Ne}, \cite{Sas}) 
is also general at the level of complexes for it automatically 
permits working with those that have quasi-coherent homology 
and not just coherent ones. Thus~$(-)^!$ is realized as  
a $\Dqc^+$-valued functor and its characterizing 
properties are all expressed in $\Dqc^+$ terms.
(Here $\Dqc^+$ has the usual meaning, namely 
the derived category of complexes $E$
such that $H^n(E)$ is quasi-coherent for all~$n$ 
and vanishes for $n\ll 0$.)

In this paper we generalize
Nagata's theorem to the category~$\Se$ of separated 
essentially finite type maps of noetherian schemes and 
in the process we also extend a few other general results 
from the finite-type setup to the essentially finite-type one.
A particular consequence, which is one of the main 
motivations behind this paper, is that for ordinary noetherian
schemes, $(-)^!$ can be defined over all of~$\Se$.  

By definition, a map in~$\Se$, 
over sufficiently small affine open subsets 
of the base and source, corresponds to 
ring homomorphisms that are essentially of finite-type, 
see \ref{def:essloc}(a).
Keeping aside considerations from duality, 
it seems reasonable to clarify at the outset, 
what the right notion of compactifiability in~$\Se$ should 
be. Specifically, let us address the question:
how should the notions of proper maps
and open immersions be generalized in~$\Se$? 
The answers are perhaps not so surprising but 
are also not so obvious at first glance.
We show below that the universally closed morphisms 
in~$\Se$ are automatically 
proper (Remark~\ref{rem:proper}). Moreover, we 
define the class of localizing immersions---these
are injective localizing morphisms, i.e., morphisms that 
are one-to-one and on sufficiently small
open sets on the base and source correspond to localization 
of rings, (see~\ref{def:locimm}, \ref{sssec:immersion})---and we
show that they give a natural generalization of open immersions.
(In particular, finite-type + localizing immersion = open immersion,
see~\ref{sssec:finitetype}).
This is convenient from the perspective of 
duality theory because the candidates for~$f^!$ in both the
cases are already available, namely, the right adjoint
$f^{\times}$ to the derived direct image 
$\R f_*$ if~$f$ is proper and the inverse image~$f^*$ 
if~$f$ is a localizing map,
and the required compatiblities between~$(-)^{\times}$
and~$(-)^*$ are generally known.  

Our ``essential" compactification theorem states 
that any map in~$\Se$ 
factors as a localizing immersion followed by a proper map 
(Theorem~\ref{cor:Nagata}). We show 
this by putting together Nagata's theorem from the finite-type case 
and the key global factorization theorem 
namely, that any $\Se$-map
factors as a localizing immersion followed by a finite-type map
in~$\Se$ (Theorem~\ref{thm:factorization}). 
Such a factorization is available locally by definition,
but it is not clear why one can find it globally.
We should also mention here that any attempt 
to mimic the known proofs of 
Nagata's compactification theorem seems difficult to carry out 
because working with localizing immersions as opposed
to open immersions has certain limitations, 
see \ref{ssec:pathos}.

Our methods also apply to produce extensions of 
other general results from the finite-type case 
to the essentially finite-type case. These include Chow's lemma,
Zariski's main theorem, elimination of indeterminacies of a rational
map via blow-ups and some results concerning birational maps,
see~\S\ref{sec:rational}.
These extensions were inspired by Conrad's exposition 
(\cite{Con2}) of Deligne's notes on Nagata's theorem.
We also prove other global factorization results which 
show that an $\Se$-map with a given property factors 
globally as a localizing immersion followed by a finite-type map 
with the same property (see~\ref{prop:properties}).

Finally, in our applications to duality, apart from the basic one,
namely that $(-)^!$ is defined over all of $\Se$
(Theorem~\ref{thm:uppershriek}), there is 
also a characterization of essentially perfect maps 
($\Se$-maps of finite tor dimension) in terms of relative
dualizing complexes, such as, $f \colon X \to Y$ is essentially
perfect iff the natural map is an isomorphism 
\smash{$f^!\OY \Otimes \bL f^*(-) \iso f^!(-)$} on~$\Dqc^+(Y)$,
see Theorem~\ref{thm:perfect}.

Prior to this paper, the most general results on duality 
over~$\Se$ would come by applying the main results 
of~\cite{Ny} or by assuming additional hypothesis that 
the underlying schemes admit dualizing complexes and 
working with~$\Dc^+$ instead of~$\Dqc^+$, 
see \cite[9.2,9.3]{DCC} for instance. 
Recent work of Yekutieli and Zhang \cite{YZ} 
on rigid dualizing complexes also develops 
duality over essentially finite-type maps. But their 
emphasis is different and the scope of their work is 
limited compared to ours. 

\section{Preliminaries on localizing immersions}
\label{sec:Prelim}

Localizing immersions are defined in \ref{def:locimm} below.
This basic notion is used throughout the paper. In \ref{ssec:props}
we give its elementary properties which show that in many respects
localizing immersions behave like open immersions. 
However, also see~\ref{ssec:pathos}. 

Let $\phi \colon A \to B$ be a ring homomorphism. We will call~$\phi$
a localizing homomorphism if~$B$ is a localization of~$A$ with 
$\phi$ as the canonical map.

\begin{adefi}
\label{def:essloc}
\hfill

(a) A map $f \colon X \to Y$ of noetherian schemes is 
said to be \emph{essentially of finite type} 
if every point $y \in Y$ has an affine open neighborhood
$V = \Spec(A)$ such that $f^{-1}V$
is covered by finitely many affine open 
$U_i = \Spec(B_i)$
for which the corresponding ring homomorphisms 
$\phi_i \colon A\to B_i$ are essentially of finite type.

(b) If, moreover, in (a), each~$\phi_i$ is a localizing 
homomorphism, then we say that $f$ is \emph{localizing}.
\end{adefi}

\sss \label{ssec:compbasech}
The defining properties in (a) and (b) behave well
with respect to composition and base change.
Thus if $f \colon X \to Y$ and $g \colon Y \to Z$ are maps 
of noetherian schemes and if both~$f$ and~$g$ are 
essentially of finite type, then so is the composition~$gf$;
if~$gf$ and~$g$ are essentially of finite type then so is~$f$;
and if $Y'\to Y$ is any scheme-map with $Y'$ noetherian, then 
$X' \set Y' \times_Y X$ is noetherian and the natural projection 
$f' \colon X' \to Y'$ is essentially of finite type. Substituting
``localizing" in place of ``essentially of finite-type" everywhere
gives similar valid statements.

\sss \label{ssec:notXlocal} 
If $f \colon X \to Y$ is essentially of finite type
with~$Y$ affine, say $Y = \Spec(A)$, then~$X$ is 
covered by finitely many affine open $V_{j} = \Spec(B_{j})$
such that the associated ring homomorphism $\phi_j \colon A\to B_{j}$
is essentially of finite type. However it is not clear and may not
be true in general that for \emph{every} affine open $V = \Spec(B)$
in~$X$, the corresponding ring homomorphism 
$\phi \colon A \to B$ is essentially of finite type. 
As before, substituting ``localizing" in place 
of ``essentially finite-type" everywhere
gives analogous statements.

\sss \label{ssec:examples}
Let us consider some examples of localizing maps.
For any scheme~$Y$ 
and any point $y \in Y$, the natural map 
$\Spec(\OYy) \to Y$ is localizing. In the 
finite-type case, open immersions are localizing maps but 
not conversely. Here are some general examples.
We will start with a base scheme~$Y$.

$\bullet$ Let $U,V$ be open subschemes of~$Y$. Let~$X$ be 
the gluing of~$U$ and~$V$ along a nonempty open 
subset of~$U \cap V$.
Then the natural map $X \to Y$ is localizing but 
not separated in general. 

$\bullet$ Let $\{U_i\}$ be a finite collection of open subsets of~$Y$.
Let $X$ be the disjoint union $\coprod_iU_i$. Then the natural 
map $X \to Y$ is localizing and separated but is not an open 
immersion in general.

In both these examples $X \to Y$ is not one-to-one. 
It is easy to show that a \emph{finite-type} localizing map that is also 
set-theoretically injective, is an open immersion, 
(\ref{sssec:finitetype}). This motivates our
choice of considering the notion of injective localizing morphisms  
as a natural extension of that of open immersions when going
from the finite-type to the essentially finite-type case. This choice
will be further vindicated by the properties and theorems that will follow.

\sss \label{ssec:genericsubsch}
Let $X$ be a noetherian scheme. The subset~$G_X$ consisting of 
all the generic points of~$X$ is finite and discrete and equals the
intersection of all the dense open subsets of~$X$.
Via restriction from~$X$ we may equip~$G_X$ with the sheaf 
of rings $\OX|_{G_X}$. Clearly,  
$G_X \cong \coprod_{\gamma \in G_X} \Spec(\cO_{X, \gamma})$.
In particular, $G_X$ forms an Artinian affine scheme. 
The natural induced scheme-map $G_X \to X$ is localizing
because it is evidently so on each component 
$\Spec(\cO_{X, \gamma})$. 

For any noetherian scheme $X$, the associated Artinian 
scheme~$G_X$ will be called as the \emph{generic subscheme} of~$X$. 
If $f \colon X \to Y$ is localizing, then the generic points 
of~$X$ map to those of~$Y$ and each point (component) of~$G_X$ 
maps isomorphically to its image in~$G_Y$. 

\begin{alem}
\label{lem:injloc}
Let $f \colon X \to Y$ be localizing. Then the following are 
equivalent.
\begin{enumerate}
\item The map $f$ is set-theoretically injective.
\item The map $f$ is separated and sends $G_X$ injectively 
inside $G_Y$.
\end{enumerate}
\end{alem}

\begin{proof}
(i) $\Lra$ (ii). We use the valuative criterion to check that $f$ is separated.
Let~$V$ be a valuation ring with quotient field~$K$. Suppose
there exist maps $\alpha \colon \Spec(K) \to X$ and
$\beta \colon \Spec(V) \to Y$ such that~$\beta$ restricts 
to~$f\alpha$ on~$\Spec(K)$. For $i = 1,2$, let 
$\alpha_i \colon \Spec(V) \to X$ be maps lifting~$\alpha$ that agree
with~$\beta$ on~$\Spec(K)$. Since~$f$ is injective, the $\alpha_i$'s
agree with each other set-theoretically. 
Let $x \in X$ be the image of the closed 
point of~$V$ under~$\alpha_i$, so that~$y = f(x)$ is the image 
of the closed point under~$\alpha$. Since~$V$ is local, $\alpha_i$
factors through the natural map $\Spec(\OXx) \to X$ while~$\alpha$
factors through $\Spec(\OYy) \to Y$. Since~$f$ is localizing,
hence the natural map $\OYy \to \OXx$ is an isomorphism, 
whence $\alpha_1 = \alpha_2$. By the valuative criterion, $f$
is separated.
  
(ii) $\Lra$ (i). Suppose $y \in Y$ has two distinct 
points, say $x_1, x_2$, in its
preimage. Since $f$ is localizing, there exist open neighborhoods
$V_i$ of $x_i$ respectively such that $f\big|_{V_i}$ is injective. 
Since $G_X$ injects into $G_Y$, hence neither
$y$ nor the $x_i$'s are generic and moreover $V \set V_1 \cap V_2$
is a nonempty open set whose image in~$Y$ contains all the 
generic points that specialize to~$y$. Since $V_1 \cup V_2$ 
is separated over $Y$ hence the natural 
map $\phi \colon V \to V_1 \times_Y V_2$ is a closed immersion.
By applying a base-change $Y' = \Spec(\OYy) \to Y$, we obtain a
closed immersion $\phi' \colon V' \to V_1' \times_Y' V_2'$ which is 
also an open immersion since $V_i' \iso Y'$. Since $Y'$ is connected,
we obtain $V' = V_i'$, thus forcing $x_1 = x_2$, a contradiction. 
 \end{proof}

\begin{adefi}
\label{def:locimm}
We define a map $f \colon X \to Y$ of noetherian schemes to be 
a localizing immersion, if $f$ satisfies the two equivalent conditions 
of Lemma~\ref{lem:injloc} above. 
\end{adefi}

\sss \label{ssec:props}
Here are some basic properties of localizing immersions. 
In what follows, all schemes shall be tacitly assumed to be noetherian.

\ssss \label{sssec:compbasech}
\emph{The property of being a localizing immersion 
behaves well under compositions and base-change. Thus the 
assertions in \textup{\ref{ssec:compbasech}} hold 
with ``localizing immersion" in place of ``essentially of finite type".}
The assertions about compositions are obvious. 
Regarding base-change, note that 
if $f \colon X \to Y$ is a localizing immersion, then 
for any $y \in Y$, the fiber-map $X_y \to \Spec(k(y))$ is an 
isomorphism. This property on fibers carries over to $f' \colon X' \to Y'$,
whence $f'$ is injective.

\ssss \label{sssec:immersion}
\emph{If $f \colon X \to Y$ is a localizing immersion, 
then~$f(X)$ is stable under generization in~$Y$ and the 
natural map of ringed-spaces induces an isomorphism
$\psi \colon (X, \OX) \iso (f(X), \OY|_{f(X)})$.} 
All the assertions can be checked locally on~$Y$, hence we 
assume that $Y$ is affine, say $\Spec(A)$. Let us write~$X$ 
as a finite union of open subschemes $U_i$, each of the form
$U_i \cong \Spec(S_i^{-1}A)$. Since each~$f(U_i)$ is stable 
under generization, so is the union~$f(X)$. The natural 
topological map 
$|\psi| \colon X \to f(X)$ is a continuous bijection 
which is a homoemorphism when restricted to each~$U_i$ 
and hence~$|\psi|$ is a homoemorphism. Since~$\psi$ is an
isomorphism over each $U_i$, it so globally.

\ssss \label{sssec:finitetype}
\emph{A localizing immersion that is of finite-type
is an open immersion.}
For a ring~$A$ and a multiplicative subset $S \subset A$,
if $S^{-1}A$ is finitely generated over~$A$, then 
there exists an element $s \in S$ such that the natural map
$A[1/s] \to S^{-1}A$ is an isomorphism.
Thus if $f \colon X \to Y$ is of finite type and 
localizing, then over sufficiently small open subsets of~$X$,
$f$ is an open immersion. Hence $f(X)$ is an open subset of~$Y$.
If~$f$ is also injective, then by~\ref{sssec:immersion}, it
is an open immersion.

\ssss \label{sssec:fiberprod}
\emph{In the situation of \textup{\ref{sssec:immersion}},
if $Y'$ is a noetherian scheme over $Y$, then 
the fiber-product $X \times_Y Y'$ maps homeomorphically
to the inverse image of~$f(X)$ in~$Y'$.} This follows easily 
from~\ref{sssec:immersion}

\ssss  \label{sssec:bijective}
\emph{A surjective (and hence bijective) localizing
immersion is an isomorphism. More generally, if $f \colon X \to Y$ is 
a localizing immersion and $f(X)$ is a closed set, then $f$ is an 
isomorphism of~$X$ onto a union of connected components 
of~$Y$.} The first assertion follows immediately from~\ref{sssec:immersion}.
For the second one it suffices to show, keeping in mind
$Z = f(X)$, that 
any closed subset $Z \subset Y$ that is stable under 
generization equals a union of connected components 
of~$Y$. Since $Z^{\textup{c}} \set Y \setminus Z$ is stable 
under specialization, it contains the closure of each of the 
generic points of~$Y$ lying in it. Every point of~$Z^{\textup{c}}$ 
lies in one of these closures because $Z^{\textup{c}}$, 
being open, is stable under generization.
Thus $Z^{\textup{c}}$ is also closed.   

\ssss \label{sssec:ideal}
\emph{If $f \colon X \to Y$ is a localizing immersion, 
then any coherent ideal~$\cI$ in~$\OX$ extends to one 
in~$\OY$, i.e., there exists 
a coherent ideal~$\J$ in~$\OY$ such that $\J\OX = \cI$.}
Indeed, in view of~\ref{sssec:immersion} let us first 
think of~$X$ as a subset of~$Y$.
Then the required ideal~$\J$ is the kernel of the composition~$\phi$ of
the natural maps $\OY \to f_*\OX \onto f_*(\OX/\cI)$, because
applying the exact restriction functor~$f^{-1}$ to~$\phi$ results in 
the composition of 
\[ \OX = f^{-1}\OY \to f^{-1}f_*(\OX/\cI) = \OX/\cI, \]
which also identifies with the canonical map $\OX \to \OX/\cI$.

\ssss \label{sssec:mono}
\emph{A localizing immersion is a localizing monomorphism
and vice-versa.} Suppose $f \colon X \to Y$ is a localizing immersion
and for $i = 1,2$ there are scheme maps $g_i \colon Z \to X$ such that
$fg_1 = fg_2$. Since $f$ is set-theoretically injective, $g_1$ and~$g_2$
agree set-theoretically. In particular, for any open subset $V \subset X$,
we have $g_1^{-1}V = g_2^{-1}V$. 
It suffices to check that as~$V$ varies over sufficiently small open 
subsets of~$X$, we have $g_1 = g_2$ on $g_i^{-1}V$.
Therefore we may assume that $X,Y$ are affine and~$f$ corresponds 
to a localizing homomorphism $A \to S^{-1}A$. 
The universal property of localization now shows that~$g_1 = g_2$.
Thus~$f$ is a monomorphism.

Conversely, suppose~$f \colon X \to Y$ is a localizing monomorphism
and there are points $x_1, x_2 \in X$, $y \in Y$ such that 
$y = f(x_i)$. Since~$f$ is localizing, the induced maps on residue
fields $k(y) \to k(x_i)$ are isomorphisms. Let~$g_i$ be the composition
of the natural maps $\Spec(k(y)) \iso \Spec(k(x_i)) \to X$. Then 
$fg_1 = fg_2$ forces $g_1 = g_2$. Thus $x_1 = x_2$.

\ssss \label{sssec:quasi-affine}
\emph{A scheme map $f \colon X \to Y$ is a localizing immersion
if and only if every $y \in f(X)$ admits an affine open neighborhood
$V = \Spec(A)$ in~$Y$ such that $U = f^{-1}V$ is affine, say 
$U = \Spec(B)$, and the corresponding ring homomorphism $A \to B$
is localizing. In particular, $f$ is quasi-affine.}
It suffices to prove the ``only if" part as the remaining assertions 
follow easily. Let $f \colon X \to Y$ be a localizing immersion. As the assertion is local on~$Y$, we may assume 
that $Y$ is affine, say $Y= \Spec(A)$, to begin with. 
Let~$y$ be a point in~$f(X)$ and~$x$ its unique pre-image.
Then~$x$ has a neighborhood of the form $U = \Spec(S^{-1}A)$. 
In view of \ref{sssec:immersion}, there exists an open subset 
$V \subset Y$ such that $f^{-1}V = U$. Upon shrinking~$V$ to a 
basic affine open neighborhood $\Spec(A[1/f])$ of~$y$ 
and~$U$ to $\Spec(S^{-1}A[1/f])$, the 
assertion follows.

\sss \label{ssec:pathos}
Some aspects of handling open immersions such as
gluing of schemes do not carry over to localizing 
immersions. Here we give an example of pathological 
behavior that illustrates this.

{\small{
Let us call a subset~$Z$
of a scheme~$Y$ a \emph{localized} subset if there exists a
localizing immersion $X \to Y$ whose set-theoretic image
equals~$Z$. By \ref{sssec:immersion}, a localizing 
immersion over~$Y$ is determined up to a 
unique isomorphism by its set-theoretic image, which is  
necessarily stable under generization. However, not every 
subset of~$Y$ stable under generization is a localized one
as the following example shows.

Let $k$ be a field. Let $A = k[T_1,T_2]$. 
Set $Y = \mathbb{A}^2 = \Spec(A)$.
Let $C$ be an  irreducible curve in~$\mathbb{A}^2$
and~$p$ a closed point on ~$C$.
Let~$Z \subset \mathbb{A}^2$ 
be the complement of the set~$C'$
consisting of all the closed points of~$C$ except~$p$. 
Clearly~$Z$ is stable under generization. But~$Z$ is 
not localized. Indeed, if it were so, then~$p$ admits 
an affine open neighborhood $U = \Spec(B)$ 
in~$Z$ where~$B$ is a localization 
of~$A$. Since~$U$ is the restriction of an open 
subset of~$\mathbb{A}^2$ to~$Z$, therefore 
$U^{\text{c}} = \mathbb{A}^2 \setminus U$ is 
the union of~$C'$ and a closed subset 
of~$\mathbb{A}^2$ away from~$p$. 
Since $U^{\text{c}}$ contains finitely many points of 
codimension one and all but one of the closed points in~$C$,
it cannot be the union of hypersurfaces in~$\mathbb{A}^2$. 
This contradicts the fact that~$B$ is a localization of~$A$.

Curiously enough, $Z$ is the union of two localized subsets 
of~$\mathbb{A}^2$ and this demonstrates the problems 
in trying to glue along localizing immersions and also of 
identifying localized subsets. Let $X_1$ be the open 
subset in~$\mathbb{A}^2$ whose complement is~$C$ 
and let $X_2 = \Spec(\cO_{\mathbb{A}^2, p})$.
Since $X_1 \cup X_2 = Z$, it means the obvious naive attempt 
at glueing the~$X_i$'s along their intersection cannot 
produce a localizing immersion into~$\mathbb{A}^2$.
}}

\section{A global factorization}

The final result of this section, 
Theorem \ref{thm:factorization},
is one of the basic main results of this paper. 

As usual, all schemes will be tacitly assumed to be noetherian.

\sss \label{ssec:sch-image}
Recall that for any scheme-map $f \colon X \to Y$ 
the \emph{schematic image} of~$f$ refers to the closed 
subscheme of~$Y$ defined by the 
kernel of the natural map $\phi \colon \OY \to f_*\OX$. We say 
that~$f$ has \emph{schematically dense image} (or~$X$ has schematically dense image in~$Y$) if~$\phi$ is injective, 
i.e, the schematic image is~$Y$. If~$X$ has 
schematically dense image in~$Y$ then it also has 
topologically dense image (i.e., $\ov{f(X)} = Y$) because 
the kernel of~$\phi$ is always supported on any open 
subset away from~$f(X)$. Thus, for any $f \colon X \to Y$,
the natural factorization \smash{$X \xto{\;g\;} \ov{X} \xto{\;i\;} Y$}, 
with~$\ov{X}$ as the schematic image of~$X$, is one 
where~$g$ is schematically dense and~$i$ is a closed
immersion. The property of having schematically dense
image is preserved under flat base change, and in particular,
base change by a localizing morphism.

\begin{alem}
\label{lem:gf-to-f}
Let $X \xto{\;f\;} Y \xto{\;g\;} Z$ be scheme-maps  
such that~$gf$ is a localizing immersion and~$g$ is separated.  
If~$f$ has schematically dense image, 
then it is a localizing immersion, while if~$f$ is proper, 
it is a closed immersion. 
\end{alem}
\begin{proof}
Consider the following diagram of natural maps where 
the square is cartesian
and~$i$ satisfies $g'i = 1_X$ and $hi = f$.
\[
\begin{CD}
X @>i>> Y' @>g'>> X \\
@. @VhVV @VVgfV \\
@. Y @>g>> Z 
\end{CD}
\]
Here $h$ is a localizing immersion and since $g'$ is separated, 
$i$ is a closed immersion. 
By~\ref{sssec:ideal} we can write 
$f = i_1h_1$ where~$i_1$ is a closed immersion 
and~$h_1$ is a localizing immersion.  
If~$f$ has schematically dense 
image, then~$i_1$ is an isomorphism 
while if~$f$ is proper then~$h_1$ is a closed immersion 
by~\ref{sssec:bijective}.
\end{proof}

In view of~\ref{sssec:immersion}, 
if a localizing immersion $f \colon X \to Y$ 
has schematically dense image, %
then we also abbreviate and say that~$f$ is \emph{schematically dense} 
(or~$X$ is schematically dense in~$Y$). Schematic 
denseness of $f$ is equivalent to every open neighborhood~$U$ 
of~$f(X)$ being schematically dense in~$Y$. 
If~$X$ is schematically dense in~$Y$, then it is also so
in every open neighborhood~$U$ in~$Y$. 

Let $f \colon X \to Y$ be a localizing immersion. If~$X$ is 
topologically dense in~$Y$, then the induced map 
of generic subschemes $G_X \to G_Y$ is an
isomorphism. In general, we can always find an open 
neighborhood~$U$ of~$f(X)$ such that~$X$ is schematically 
dense in~$U$. Indeed, if $\cI$ is the kernel of the natural map 
$\OY \to f_*\OX$, then taking~$U$ to be the complement of the 
support of~$\cI$ works. For any such $U$, it holds that
$G_X \iso G_U \subset G_Y$. 

\begin{aprop}
\label{prop:approximate}
Let 
\[
\begin{CD}
X @>{f_1}>> Y_1 \\
@| @VVgV \\
X @>{f_2}>> Y_2
\end{CD}
\]
be a commutative diagram of separated scheme-maps
where each~$f_i$ is a localizing immersion and~$g$ is 
essentially of finite type. 
Then there exists an open subscheme 
$U \subset Y_1$ containing~$f_1(X)$ 
such that $g |_U \colon U \to Y_2$ is a localizing immersion.
In particular, if~$g$ is of finite type, then $g |_U$ is an open
immersion.
\end{aprop}
\begin{proof}
The last statement follows from~\ref{sssec:finitetype}. 
Replacing~$Y_i$ by open subsets~$V_i$ containing~$f_i(X)$ and 
satisfying $g(V_1) \subset V_2$ does not affect the assertion in the 
proposition. Hence by shrinking~$Y_i$'s if necessary, 
we may assume without loss of generality that~$X$ is 
schematically dense in both, $Y_1$ and~$Y_2$. 
Thus we have $G_{Y_1} \cong G_{X} \cong G_{Y_2}$. 
By Lemma~\ref{lem:injloc}, %
it suffices to find~$U$ 
containing~$f_1(X)$ such that~$g|_U$
is localizing. This problem can verified locally on~$Y_2$,
$Y_1$ and~$X$ and so we shall now assume that all three are affine.
Specifically, let $X = \Spec(B)$, $Y_i = \Spec(A_i)$, assume that 
for each~$i$, there is an isomorphism $B \cong S_i^{-1}A_i$ for a
suitable multiplicative subset~$S_i$ in~$A_i$ and assume 
that~$A_1$ is essentially of finite type over~$A_2$.

Since $f_i$ is schematically dense, 
the natural map $A_i \to B$ is injective. For convenience, 
let us identify each~$A_i$ with its image in~$B$. Thus there is a 
containment of rings $A_2 \subset A_1' \subset A_1 \subset B$ 
where $A_1'$ is finitely generated over~$A_2$ and
$A_1 = T^{-1}A_1'$ for some multiplicative subset 
$T \subset A_1'$. Let $a_1, \ldots, a_n$ be in~$A_1'$
and $x_1, \ldots, x_n$ in~$S_2$ such that 
$A_1' = A_2[a_1/x_1, \ldots, a_n/x_n]$. Then for $x = x_1 \cdots x_n$,
it holds that $A_1'[1/x] = A_2[1/x]$. Thus $A_1[1/x]$ is a localization
of~$A_2[1/x]$, so that $U = \Spec(A_1[1/x])$ 
gives the desired open set.   
\end{proof}

\begin{aprop} 
\label{prop:dominating}
Consider the following commutative diagram of separated 
scheme-maps
\[
\begin{CD}
X @>{f_1}>> Y_1 \\
@V{f_2}VV @VV{\pi_1}V \\
Y_2 @>{\pi_2}>> S
\end{CD}
\]
where~$\pi_i$ are of finite type and~$f_i$ are localizing immersions.
Then there exists a localizing immersion $h \colon X \to W$
for some scheme~$W$ of finite type over~$S$ such that~$W$ 
admits open $S$-immersions $g_i \colon W \to Y_i$
satisfying $g_ih = f_i$.
\end{aprop}

\begin{proof}
Set $P \set Y_1 \times_S Y_2$. Let~$Z$ be the schematic image
of~$X$ in~$P$ so that the obvious natural map from~$X$ to~$P$ 
factors  as \smash{$X \xto{\;j\;} Z \xto{\;k\;} P$} where~$k$  is a closed 
immersion and~$j$ has schematically dense image.
Set $g_i \set p_ik \colon Z \to Y_i$ where~$p_i$ is
the canonical projection $P \to Y_i$. 
Since $f_i = g_ij$, is a localizing immersion and~$j$ has
schematically dense image, 
hence by Lemma~\ref{lem:gf-to-f}, $j$ is a localizing immersion.
By Proposition~\ref{prop:approximate} there are open 
subschemes $U_i$ of~$Z$ containing 
$j(X)$ such that $g_i|_{U_i}$ is an open immersion. Then
$W \set U_1 \cap U_2$ gives the desired scheme satisfying the
proposition. 
\end{proof}

\begin{aprop}
\label{prop:two-glue}
Let $f \colon X \to S$ be a separated scheme map.
Suppose there are two open subsets 
$U_1, U_2$ covering~$X$ such that for each~$i$,
$f|_{U_i}$ factors as \smash{$U_i \xto{\:k_i\:} Y_i \xto{\:p_i\:} S$} 
where~$k_i$ is a localizing immersion
and~$p_i$ is separated and of finite type.
Then~$f$ admits a factorization
\smash{$X \xto{\;k\;} Y \xto{\;p\;} S$} where~$k$ is a localizing immersion
and~$p$ is separated and of finite type.
\end{aprop}

\begin{proof}
Set $U_{12} \set U_1 \cap U_2$. By Proposition \ref{prop:dominating}, 
there is scheme $Y_{12}$, a localizing immersion
$h \colon U_{12} \to Y_{12}$ and 
open immersions $g_i \colon Y_{12} \to Y_i$ such that 
$g_ih = k_i|_{U_i}$. Since each $U_i$ is homeomorphic to its 
image in~$Y_i$ (\ref{sssec:immersion}), there are
 open subsets $W_i \subset Y_i$
such that $k_i^{-1}W_i = U_{12}$. Replacing $Y_{12}$
by $g_1^{-1}W_1 \cap g_2^{-1}W_2$, we may assume that 
$Y_{12} \times_{Y_i} U_i = U_{12}$ for each~$i$. To summarize,
we now have the following commutative diagram where the two 
parallelograms on the left are cartesian, $h,k_i$ are localizing
immersions, $g_i$ are open immersions and $p_i$ are of finite type.
\[
\xymatrix{
& U_1 \ar[rr]^{k_1} & & Y_1 \ar[rd]^{p_1} \\
U_{12} \ar[ru] \ar[rd] \ar[rr]^{h} & & 
Y_{12} \ar[ru]^{g_1} \ar[rd]_{g_2} & & S  \\
& U_2 \ar[rr]_{k_2} && Y_2 \ar[ru]_{p_2}
}
\] 
The above-mentioned properties of this diagram are not affected by
any further shrinking of~$Y_{12}$ to an open subset~$Y_{12}'$ 
containing~$h(U_{12})$ and replacing~$g_i$ by $g_i|_{Y_{12}'}$.
Likewise, shrinking~$Y_1$ to an open subset~$Y_1'$,  
$Y_{12}$ to~$g_1^{-1}Y_1'$ and correspondingly 
modifying $p_1,g_1$ has no effect. A similar statement holds
for shrinking~$Y_2$. 

Since $X$ is separated over $S$, the natural immersion 
$U_{12} \to U_1 \times_S U_2$ has a closed image. However
$Y_{12} \to Y_1 \times_S Y_2$ need not be a closed immersion
and so the gluing of~$Y_i$ along~$Y_{12}$ need not be separated 
over~$S$. Our aim now is to shrink~$Y_i$'s and~$Y_{12}$ 
suitably so that separation is achieved. 

Consider the following commutative diagram where
all the squares are cartesian and where $ba$ in the bottom
row is the obvious natural map, displayed 
as factoring through the schematic 
image~$\ov{U}_{12}$ of~$U_{12}$ in $U_1 \times_S Y_2$.
\[
\begin{CD}
U_{12}' @>a'>> \ov{U}_{12}' @>b'>> U_1 \times_S U_2 
@> {\textup{projection}} >> U_2 \\
@VV{k_2'}V @VV{\ov{k}_2}V @VVV @VV{k_2}V \\
U_{12} @>a>>  \ov{U}_{12} @>b>> U_1 \times_S Y_2 
@> {\textup{projection}} >> Y_2
\end{CD}
\]
Thus $b,b'$ are closed immersions, while the vertical arrows are all
localizing immersions. Now note that~$ba$ is an immersion, 
so that~$a$ is an open immersion. 
Indeed, the natural immersion 
$U_{12} \subset U_1 \to U_1 \times_S Y_1$
also factors naturally as
$U_{12} \xto{\;e\;} U_1 \times_S Y_{12} \xto{\textup{open}} 
U_1 \times_S Y_1$
whence~$e$ is an immersion.  As 
$U_1 \times_S Y_{12} $ is open inside $U_1 \times_S Y_2$ 
therefore~$ba$ is an immersion. By definition
of~$\ov{U}_{12}$ it follows that~$a$ is a schematically dense 
open immersion. Thus~$a'$ too is a schematically dense 
open immersion.

In view of \ref{sssec:fiberprod} and \ref{sssec:bijective} we see
that~$k_2'$ is an isomorphism. Separatedness of $X/S$ implies
that $U_{12}' \cong U_{12}$ is closed in~$U_1 \times_S U_2 $,
whence~$a'$ is a closed immersion. Since~$a'$ is schematically 
dense, 
it is an isomorphism. In view of~\ref{sssec:fiberprod}
we see that the natural projection of %
$b(\ov{U}_{12} \setminus U_{12}) \subset
U_1 \times_S Y_2$ to $Y_2$ is a set $Z$ 
disjoint from~$k_2(U_2)$. As $k_2(U_2)$ is stable under 
generization (\ref{sssec:immersion}), it is also disjoint from the 
closure~$\ov{Z}$ of~$Z$ in~$Y_2$. Therefore, performing 
an open base-change by $Y_2' \set Y_2 \setminus \ov{Z} \to Y_2$ 
to the bottom row of the above diagram transforms~$a$
into an isomorphism. Let us replace~$Y_2$ by~$Y_2'$,
$Y_{12}$ by $g_2^{-1}Y_2'$, so that we may assume 
henceforth that~$a$ is an isomorphism, i.e., $U_{12}$ is a 
closed subscheme of $U_1 \times_S Y_2$.

Now consider the following diagram where too the 
squares are cartesian and $dc$ is the obvious natural immersion
displayed as factoring through the schematic 
closure~$\ov{Y}_{12}$ of~$Y_{12}$ in $Y_1 \times_S Y_2$. 
\[
\begin{CD}
Y_{12}' @>c'>> \ov{Y}_{12}' @>d'>> U_1 \times_S Y_2 
@> {\textup{projection}} >> U_1\\
@VV{k_1'}V @VV{\ov{k}_1}V @VVV @VV{k_1}V \\
Y_{12} @>c>>  \ov{Y}_{12} @>d>> Y_1 \times_S Y_2 
@>{\textup{projection}}>> Y_1
\end{CD}
\]
Thus the vertical arrows are localizing immersions, $c,c'$ are 
schematically dense open immersions, while $d,d'$ are closed immersions. Since $U_{12}$ is isomorphic to 
$Y_{12} \times_{Y_1} U_1$,
we may assume that $Y'_{12} = U_{12}$ and $k_1' = h$. 
Moreover we may
identify $c', d'$ with $a,b$ of the earlier diagram respectively.
Since~$a$ has been arranged to be an isomorphism, 
therefore~$c'$ is an isomorphism.

Arguing as before, we see that the projection of 
$d(\ov{Y}_{12} \setminus Y_{12}) \subset Y_1 \times_S Y_2$
to~$Y_1$ is a set~$W$ disjoint from~$k_1(U_1)$ and hence 
the closure~$\ov{W}$ is also disjoint from~$k_1(U_1)$.
Thus replacing~$Y_1$ with the open subset 
$Y_1 \setminus \ov{W}$ and correspondingly modifying~$Y_{12}$
ensures that~$c$ is an isomorphism. 

Since we have found suitable $Y_i$ and $Y_{12}$ for which 
$Y_{12} \to Y_1 \times_S Y_2$ is a closed immersion,
the gluing of~$Y_1$ and $Y_2$ along $Y_{12}$ yields
a scheme~$Y$ that is separated and of finite-type over~$S$. 
The natural induced map $k \colon X \to Y$ is localizing
and injective by our choice of~$Y_{12}$ whence it is a
localizing immersion as required.
\end{proof}

\begin{athm}
\label{thm:factorization}
Let $f \colon X \to S$ be a separated essentially finite-type map 
of noetherian schemes. Then~$f$ factors as 
\smash{$X \xto{\;k\;} Y \xto{\;p\;} S$} where~$k$ 
is a localizing immersion and~$p$ is separated and of finite type.
\end{athm}

\begin{proof}
By definition, there exists a finite open cover $\{U_i\}$ of~$X$  
such that $f|_{U_i}$ factors as 
\smash{$U_i \xto{\:k_i\:} Y_i \xto{\:p_i\:} S$} where~$k_i$ 
is a localizing immersion and~$p_i$ is separated and of finite type.
By Proposition~\ref{prop:two-glue} and induction 
on the number of elements in the open cover, the theorem follows. 
\end{proof}

\section{General Applications}
\label{sec:rational}

Below, we give various applications of the results from the previous 
section.

As a consequence of Theorem~\ref{thm:factorization}, 
we can now extend Nagata's compactification
theorem \cite{Ng}, \cite{Con2}, \cite{Lu} to essentially finite-type maps. 

\begin{athm}
\label{cor:Nagata}
Let $f \colon X \to S$ be a separated essentially finite-type map 
of noetherian schemes. Then~$f$ factors as 
\smash{$X \xto{\;k\;} Y \xto{\;p\;} S$} where~$k$ 
is a localizing immersion and~$p$ is proper.
\end{athm}

\begin{arem}
\label{rem:proper}
In the context of the above theorem we note that 
replacing ``finite type" with ``essentially finite type"
in the definition of properness does not define a new condition, 
i.e., a separated essentially-finite-type map $f \colon X \to S$
that is also universally closed, is necessarily proper. 
Indeed, by Theorem~\ref{thm:factorization}, 
$f$ factors as \smash{$X \xto{\;i\;} Y \xto{\;k\;} S$}
where~$i$ is a localizing immersion and~$k$ is separated and 
of finite type. Separatedness of~$k$ implies that~$i$ is universally 
closed and hence by \ref{sssec:bijective}, $i$ is a closed immersion. 
Thus~$f$ is of finite type and hence proper.
\end{arem}

The remaining general applications in this section 
will play no role in the results of \S~\ref{sec:duality}.

Next let us look at Zariski's Main Theorem. 
We define a map $f \colon X \to Y$ of noetherian
schemes to be \emph{essentially quasi-finite}, if it is essentially 
of finite type and for any $y \in Y$, the fiber~$X_y$ is 
algebraically finite, i.e., $X_y \cong \Spec(A)$ where~$A$ 
has finite vector-space dimension over the 
residue field~$k(y)$ at~$y$.

\begin{athm}
\label{thm:ZMT}
Let $f \colon X \to Y$ be a separated essentially quasi-finite 
map of noetherian schemes. Then~$f$ factors as 
\smash{$X \xto{\; i\;} Z \xto{\; h\;} Y$} where~$i$ is a localizing 
immersion and~$g$ is finite.
\end{athm}

\begin{proof}
By Theorem~\ref{thm:factorization}, 
we can factor $f$ as $X \xto{\; j \;} Z \xto{\; g \;} Y$
where~$j$ is a localizing immersion and~$g$ is separated 
and of finite type.
Pick a point $x \in X$ and let $z = j(x)$, $y = f(x)$. Let $G = Z_y$
and $F = X_y$ be the corresponding fibers over $y$. 
Since~$G$ is of finite type over~$k(y)$ and 
$\cO_{G,z} \cong \cO_{F,x}$ is a finite $k(y)$-module, 
therefore~$z$ is an isolated point of~$G$. Thus the set~$U$ of 
points in~$Z$ that are isolated in their fiber over~$Y$ 
contains~$j(X)$ and by \cite[13.1.4]{EGAIV}, 
$U$ is open. Replacing~$Z$ by~$U$
(and~$g$ by~$g|_U$) we may therefore assume 
that the fibers of~$g$ are discrete 
and hence that~$g$ is quasi-finite. 
By Zariski's Main Theorem, $g$ factors 
as \smash{$Z \xto{\;l\;} \ov{Z} \xto{\;h\;} Y$} where~$l$ is 
open and~$h$ finite. Using~$\ov{Z}$ in place of~$Z$ and 
setting $i = lj$ we deduce the theorem.
\end{proof}

The rest of this section concerns some general results on 
rational maps and blow-ups. These generalize statements 
in the finite-type case,
all of which can all be found in \cite{Con2}. Let us 
fix some terminology and recall some general facts about 
blow-ups.

\sss
\label{ssec:blowuppreview}
Let $f \colon X \to Y$ be a localizing immersion
and $\cI$ a coherent ideal in~$\OY$. 
Let $\wt{Y} = \Bl_{\cI}(Y)$ be the blow-up of~$Y$
with respect to~$\cI$. We say that~$\wt{Y}$ is an 
\emph{$X$-admissible} blow-up of~$Y$ (or the blow-up map
$\pi \colon \wt{Y} \to Y$ is $X$-admissible) if the closed 
subscheme of~$Y$ defined by~$\cI$ is
disjoint from~$f(X)$. This is equivalent to the condition 
that for any $y \in f(X)$, we have $\cI_y = \OYy$ and is 
also equivalent to requiring that there 
exist an open neighborhood~$U$ of~$f(X)$ such that 
$\cI|_U = \OU$. In such a situation, $f$ lifts to the
blow-up~$\wt{Y}$, i.e., there is a localizing immersion
$\wt{f} \colon X \to \wt{Y}$ such that $\pi \wt{f} = f$. 

%
%

For a scheme map $f \colon X \to Y$ and blow-ups
$\pi_1 \colon \wt{X} \to X$, $\pi_2 \colon \wt{Y} \to Y$, we say that a map 
$f' \colon \wt{X} \to \wt{Y}$ lifts~$f$ if the following 
diagram commutes. 
\[
\begin{CD}
\wt{X} @>{f'}>> \wt{Y} \\
@V{\pi_1}VV @VV{\pi_2}V \\
X @>f>> Y
\end{CD}
\]
If we fix the blow-up of~$Y$, say $\wt{Y} = \Bl_{\cI}(Y)$, then 
any lift~$f'$ of~$f$ factors through the canonical map 
$\wt{f} \colon \Bl_{\cI\OX}(X) \to \Bl_{\cI}(Y)$.
If $f$ is a localizing immersion, 
then with $f'=\wt{f}$, the above diagram is cartesian as $f$ is flat. 
In this case, if~$\pi_2$ is $X$-admissible, 
then~$\pi_1$ is an isomorphism, so our
notation of~$\wt{f}$ remains consistent with that of the previous 
paragraph. 

If~$U$ is an open
subscheme of a noetherian scheme~$X$ then~$U$ is schematically
dense in $\wt{X} = \Bl_{\I}(X)$ where $\I$ is any coherent ideal defining
the closed set $X \setminus U$. This follows from the fact that 
the complement of~$U$ in~$\wt{X}$ is defined by an invertible ideal,
namely~$\wt{\I}$. If~$U$ is dense (resp. \! schematically dense) in~$X$,
then it remains so in any $U$-admissible blow-up 
$\wt{X} = \Bl_{\cI}(X)$ as can be seen by looking at the open immersions
$U \to X \setminus V(\cI) = \wt{X} \setminus V(\cI) \to \wt{X}$,
all of which are dense (resp. \! schematically dense). 

The statements in the previous paragraph can be generalized to 
localizing immersions because if $f \colon Z \to X$ is a 
localizing immersion, then there is an open subset $U \subset X$ 
containing~$f(Z)$ such that~$Z$ is schematically dense in~$U$.
Thus there is a blow-up $\wt{X}$ such that $Z$ 
is schematically dense in~$\wt{X}$ namely, 
the blow-up of $X$ along $X \setminus U$. Also, if~$Z$
is dense (or schematically dense) in~$X$, then it remains so
in any $Z$-admissible blow-up of~$X$.

A composition of blow-ups is again a blow-up. 
In fact, for an open inclusion $U \subset X$, 
let \smash{$X'' \xto{\;q\;} X' \xto{\;p\;} X$} be 
$U$-admissible blow-ups (for~$q$, note that
$p^{-1}U \cong U$). 
Then~$pq$ is also a blow-up, which is necessarily 
$U$-admissible (see \cite[Lemma 1.2]{Con2}).

The following theorem generalizes the result that a ``birational" map
can be transformed into an open immersion after suitably blowing 
up the base and source, (\cite[Corollary 4.4]{Con2}).
The finite-type version there actually seems to follow 
immediately from [ibid, Lemma 2.7] itself and also 
appears in \cite[Proposition 3, p.\ 30]{Ray}.

\begin{athm}
\label{thm:birational}
Let $X \xto{\;f\;} Y \xto{\;g\;} Z$ be separated essentially finite-type
maps of noetherian schemes such that~$f$ and~$gf$ are 
localizing immersions. Then there are $X$-admissible blow-ups
$\wt{Y}$, $\wt{Z}$ and a lift $g' \colon \wt{Y} \to \wt{Z}$ of~$g$
such that~$g'$ is a localizing immersion. 
\end{athm}

\begin{arem}
\label{rem:properbirational}
(The proper case) Assume that $X$ is dense in $Z$. If~$g$ is proper,
then so is~$g'$ so that by \ref{sssec:bijective}, $g'$ is an isomorphism.
\end{arem}

\begin{proof}
By Theorem~\ref{thm:factorization}, $g$ factors as 
$Y \xto{\;k\;} Y_1 \xto{\;p\;} Z$ 
where~$k$ is a localizing immersion and $p$ is separated and 
of finite type. By Proposition~\ref{prop:dominating}, $p$ is a 
birational map, i.e., there exists an 
open subset $U \subset Y_1$
containing $kf(X)$ such that~$p|_U$ is an open immersion. 
If~$U$ is dense in~$Y_1$ and~$Z$, then by 
\cite[Corollary 4.4]{Con2} there is a $U$-admissible 
blow-up of~$Z$ such that 
the natural induced map $\wt{p} \colon \wt{Y}_1 \to \wt{Z}$
is an open immersion. In the general case, keeping in mind that 
a composition of $U$-admissible blow-ups is again a
($U$-admissible) blow-up, we first blow-up 
$Y_1$ and~$Z$ along the largest coherent ideal defining
$Y_1 \setminus U$ and $Z \setminus U$ respectively,
so that $U$ becomes schematically dense in $Y_1, Z$ and then 
use \cite[Corollary 4.4]{Con2}. Thus we have an open immersion 
$p' \colon \wt{Y}_1 \to \wt{Z}$ that lifts~$p$.
The $U$-admissible blow-up $\wt{Y}_1$ induces an
$X$-admissible blow-up~$\wt{Y}$of~$Y$ and composing the 
natural map $\wt{k} \colon \wt{Y} \to \wt{Y}_1$ with~$p'$
gives us a lifting~$g'$ of~$g$ as desired. 
\end{proof}

Next we consider the gluing of two or more schemes 
(of essentially finite-type over a base scheme~$S$)
along open immersions coming from a common scheme. 
Such a gluing need not result in a scheme separated over 
the base but after suitable $U$-admissible blow-ups it does.
The case of gluing two schemes 
can be thought of as an abstract version of Chow's lemma.

\begin{aprop}
\label{prop:glueblowups}
Let $S$ be a noetherian scheme and let $f_i \colon U \to Y_i$
be a finite collection of localizing $S$-immersions of essentially 
finite type $S$-schemes. 
Then there exist $U$-admissible blow-ups 
$\wt{Y_i} \to Y_i$, a separated finite-type $S$-scheme~$Z$ 
and dense localizing $S$-immersions 
$g_i \colon \wt{Y_i} \to Z$
such that for all $i,j$ we have 
$g_i\wt{f}_i = g_j\wt{f}_j$.
\end{aprop}

\begin{proof}
It suffices to prove the result in the case where there are only 
two~$Y_i$'s since the general case follows by induction in view of
the fact that composition of $U$-admissible blow-ups is again a
$U$-admissible blow-up.
By Theorem~\ref{thm:factorization}, 
the natural map $Y_i \to S$ factors as 
\smash{$Y_i \xto{k_i} Z_i \xto{p_i} S$}
where~$k_i$ is a localizing immersion and~$p_i$ is separated and 
of finite type. By Proposition~\ref{prop:dominating},
there exists a localizing immersion
$h \colon X \to W$ for some scheme~$W$ of
finite type over~$S$ and there are open $S$-immersions
$e_i \colon W \to Z_i$ such that $e_ih = g_if_i$. 
By \cite[Corollary 2.10]{Con2} there 
exist $W$-admissible blow-ups $\wt{Z_i} \to Z_i$ and
a separated finite-type $S$-scheme~$Z$ together 
with open immersions $\wt{Z_i} \to Z$ that agree on~$W$. 
By construction, these open immersions can also be arranged to 
be dense.
The blow-ups on~$Z_i$ induce corresponding ones on~$Y_i$
so setting~$g_i$ to be the composition of the natural maps
$\wt{Y}_i \to \wt{Z}_i \to Z$ proves the theorem.
\end{proof}

Let us generalize Chow's Lemma. We call a map $Z \to S$ of 
noetherian schemes \emph{essentially quasi-projective} if it 
factors as \smash{$Z \xto{\;j\;} Y \xto{\;p\;} S$} 
where~$j$ is a localizing immersion
and~$p$ is a projective morphism.  Here we use 
the definition of quasi-projectivity and
projectivity as in \cite[\S 5.3, 5.5]{EGAII} 

\begin{athm}
\label{thm:Chow}
Let $\pi \colon Y \to S$ be a separated essentially-finite-type morphism 
of noetherian schemes. Let $f \colon U \to Y$ be a dense localizing  
immersion such that the natural map $\pi f \colon U \to S$ 
is essentially quasi-projective. Then there exists an $U$-admissible
blow-up $\wt{Y} \to Y$ such that the natural map~$\wt{Y} \to S$ 
is essentially quasi-projective.
\end{athm}

\begin{proof}
Let $U \xto{\;f'\;} U^* \xto{\;\pi'\;} S$ be a factorization of~$\pi f$
with~$f'$ a localizing immersion and 
$\pi'$ a projective morphism.
By Proposition~\ref{prop:glueblowups}, with $Y_1 = Y$, 
$Y_2 = U^*$ and $f_1 = f$, $f_2 = f'$, we find that for
suitable $U$-admissible blow-ups of~$Y$ and~$U^*$, 
there are dense localizing immersions
$g_1 \colon \wt{Y} \to Z$, $g_2 \colon \wt{U}^* \to Z$ 
that agree on~$U$ where~$Z$ is separated 
and of finite-type over~$S$. Since~$g_2$ is dense and proper, by \ref{sssec:bijective}, it is an isomorphism. Thus~$Z$ is projective
over~$S$ and~$\wt{Y}$ essentially quasi-projective.
\end{proof}

The finite-type version of Theorem~\ref{thm:birational}  
is closely related to a general result (\cite[Theorem 2.4]{Con2}) about
eliminating indeterminacy of rational maps via blow-ups.
This result is one of the important steps in the proof of 
Nagata's compactification theorem.
Here is the analogous result for essentially finite type maps.

\begin{athm}
\label{thm:quasi-dom}
Let $S$ be a noetherian scheme and let $X,Y$ be noetherian 
schemes separated and of essentially finite type over~$S$. 
Suppose there are $S$-morphisms 
\smash{$X \xgets{j} Z \xto{f} Y$} with~$j$
a localizing immersion. Then there 
exists a $Z$-admissible blow-up $\wt{X}$, a localizing
immersion $j' \colon Z' \to \wt{X}$ extending the natural inclusion
$\wt{j} \colon Z \to \wt{X}$ and an $S$-morphism
$f' \colon Z' \to Y$ extending~$f$, such that the natural map 
$Z' \to \wt{X} \times_S Y$ is a closed immersion.
\end{athm}

\begin{arem}
\label{rem:indeterminacy}
Suppose $Y \to S$ is proper and~$Z$ is 
dense in~$X$ so that~$Z$, $Z'$ are dense 
in~$\wt{X}$. 
Properness of the composition 
$Z' \to \wt{X} \times_S Y \to X$ implies that $Z' = \wt{X}$. Thus the 
domain of the rational map~$f$ extends to the whole of~$\wt{X}$ 
in this case. 
\end{arem}

One way of proving Theorem~\ref{thm:quasi-dom} would be 
to reduce it to the original result from the finite-type case
via our main factorization theorem
(\ref{thm:factorization}). Though this can be carried out,
we use a somewhat different approach below and deduce it using 
Theorem~\ref{thm:birational}. This illustrates the close relation
between these results.

\begin{proof}
For any choice of $Z', \wt{X}$, specifying an $S$-map $Z' \to Y$
is equivalent to specifying a map $Z' \to X \times_S Y$. Thus 
replacing~$S$ by~$X$ and~$Y$ by~$X \times_S Y$ does not 
change the problem. Hence we shall from now on assume 
that $S = X$. In particular, there is now a natural map 
$g \colon Y \to X$.

First we find a $Z$-admissible blow-up $\wt{X} \to X$ such 
that~$Z$ is schematically dense in~$\wt{X}$. 
Next note that replacing~$X$ 
by~$\wt{X}$ and~$Y$ by $\wt{X} \times_X Y$ does not 
affect the problem. Hence we can and will now assume 
that $Z$ is schematically dense in~$X$. This assumption 
is not affected by any further $Z$-admissible blow-up of~$X$.
Therefore for any choice of a further blow-up~$\wt{X}$ and 
any possible~$f'$ extending~$f$, the schematic image 
of~$f'$ equals that of~$f$. Hence, by replacing~$Y$ 
with the schematic image of~$f$, we may assume without
loss of generality that $f$ has schematically dense image.  

Since $j=gf$ is a localizing immersion, by Lemma~\ref{lem:gf-to-f}, 
$f$ is a localizing immersion. Hence by 
Theorem~\ref{thm:birational}, there exists a $Z$-admissible 
blow-up of~$X$ such that the natural induced 
map $\wt{g} \colon \wt{Y} \to \wt{X}$ is a localizing immersion.
We choose $Z' = \wt{Y}$. It remains to verify that the natural
map $\wt{Y} \xto{\;h\;} \wt{X} \times_X Y$ is a closed immersion.
Let $p$ denote the canonical projection $\wt{X} \times_X Y \to \wt{X}$.
Since~$h$ is proper and $ph = \wt{g}$ is a localizing immersion, 
by Lemma~\ref{lem:gf-to-f}, $h$ is a closed immersion.
\end{proof}

\section{Duality for essentially-finite-type maps}
\label{sec:duality}

The main applications to duality given here, Theorems 
\ref{thm:uppershriek} and \ref{thm:perfect}, follow quickly from
the existence of essential compactifications 
(Theorem~\ref{cor:Nagata}) 
since the proofs are already there in literature.
Towards the end we show that many properties of 
essentially finite-type maps can be approximated 
by finite-type ones by means of a global
factorization result, see Proposition~\ref{prop:properties}.

\sss
\label{ssubec:derivedcat}
Let us recall some basic notation. For a scheme~$X$, we 
use $\D(X)$ to denote the derived category of the 
category of~$\OX$-modules, $\Dqc(X)$ (resp. \!$\Dc(X)$) to
denote the full subcategory whose objects are the complexes
having quasi-coherent
(resp. \!coherent) homology and $\D^+_{*}(X)$ 
(resp. \!$\D^-_{*}(X)$, resp. \!$\D^{\textup{b}}_{*}(X)$) to denote 
the full subcategory of $\D_{*}(X)$ whose objects are 
complexes $\F$ such that $H^n\F = 0$ for $n\gg0$ 
(resp. \!$n\ll0$, resp. \!$|n|\gg0$).

Recall that $\Se$ is the category of essentially finite-type 
morphisms of noetherian schemes.
An \emph{essentially \'{e}tale} map of noetherian schemes is 
a separated formally \'{e}tale map that is essentially of finite type.
Essentially \'{e}tale maps form a larger subcategory of~$\Se$ 
than the one of localizing maps.

\sss
\label{ssec:pastepseudofunctors}
One of the main results in \cite[Theorem 7.1.6]{Ny} 
was that for ordinary (noetherian) 
schemes, a pseudofunctorial construction of $(-)^!$ satisfying
a flat-base-change isomorphism is valid over the category~$C$
of composites of proper maps and essentially \'{e}tale 
maps.\footnote{The statement given there involves only finite-type
\'{e}tale maps, but the proof works under the 
essentially finite-type hypothesis too.} 
In particular, the theorem shows that the existence of essential compactifcations is not needed for pseudofunctorially defining 
$(-)^!$ over~$C \subset \Se$, but it does not prove
that~$(-)^!$ is defined over all of~$\Se$. 
Theorem~\ref{cor:Nagata} completes the picture 
because now we know that $C = \Se$. 
For convenience, let us briefly recall 
the basic defining properties of~$(-)^!$.

$(-)^!$ is a
contravariant $\Dqc^+$-valued pseudofunctor on~$C$ such that:
\begin{itemize}
\item on proper maps, $(-)^!$ is pseudofunctorially 
isomorphic to the right adjoint of the right-derived direct image 
pseudofunctor $\R f_*$;
\item on essentially \'{e}tale maps, $(-)^!$ equals the inverse image 
pseudofunctor $(-)^*$;
\item for any fibered square $\sfr$ of morphisms of noetherian schemes as follows 
\[
\begin{CD}
U @>{j}>> X \\
@V{g}VV @VV{f}V \\
V @>{i}>> Y
\end{CD}
\]
where $f$ is proper (and hence is in~$\Se$) 
and $i$ is flat, there is a 
flat-base-change isomorphism $\beta_{\sfr} \colon j^*f^! \iso g^!i^*$,
(see \cite[4.4.3]{Lip-LNM}). Moreover, if $i,j$ are essentially
\'{e}tale, then $\beta_{\sfr}$ agrees with the natural isomorphisms
\[
j^*f^! \iso j^!f^! \iso (fj)^! = (ig)^! \iso g^!i^! \iso g^!i^*.
\]
\end{itemize}

These properties together with the existence of essential compactifications extend the abstract theory of $(-)^!$ to all
of~$\Se$:

\begin{athm}
\label{thm:uppershriek}
The Grothendieck duality pseudofunctor $(-)^!$ 
exists on the entire category $\Se$ of separated 
essentially-finite-type maps of noetherian schemes
and satisfies compatibility with flat base change.
\end{athm}

\sss
\label{ssec:Verdiercoherence}
For smooth or finite maps, $(-)^!$ has a
concrete descriptions. We only describe the 
functorial aspects and not the pseudofunctorial one 
over these subcategories. 
 
Let $f \colon X \to Y$ be a closed immersion.
Since $f$ is proper, $f^! \cong f^{\times}$
the right adjoint of $\R f_*$, so 
one can use $f^!(-) \cong f^{-1}\R\sHom_{\OY}(f_*\OX, -)$, 
(see for instance, \cite[p.\ 172, 6.8]{RD}). 
In particular, if $f$ is a regular immersion induced locally
by regular sequences of length~$d$,
then we also have the fundamental local isomorphism
\smash{$f^!(-) \cong  \omega[-d] \Otimes \bL f^*(-)$} 
([ibid. p.\ 180, 7.3]), %
where $\omega$ is the $d$-th exterior power of the 
normal module for~$f$. (Also see \cite[p.~111]{Lip-Ast} for a 
more intrinsic description of this isomorphism.)

Recall that a scheme-map is called \emph{essentially smooth} 
if it is in~$\Se$ and it is formally smooth. Any such map, 
say~$f$, has a locally constant relative dimension that 
corresponds to the rank of the module of relative 
differentials~$\Omega_f$, (\cite[16.10.2]{EGAIV}).
If $f \colon X \to Y$ is essentially smooth of 
relative dimension~$d$, then Verdier's argument 
in \cite[Theorem 3]{V}
shows that there is a natural 
isomorphism $f^!(-) \cong f^*(-) \otimes_X \Omega^d[d]$
and in particular that $f^! \OY \cong \Omega^d[d]$. 
(Verdier's argument can be carried out quite formally and the
main input is the fundamental local isomorphism above
applied to the diagonal map $X \to X \times_Y X$ which 
is locally given by a regular sequence of length~$d$.)

It follows that for any $\Se$-map $f \colon X \to Y$, %
we have $f^!\Dc^+(Y) \subset \Dc^+(X)$. Indeed, since $f$
factors locally on~$X$ as a closed immersion into a smooth 
map, and in these cases the above formulas for~$f^!$ 
preserve coherence of homology, the same holds for 
the general case.

\sss
\label{ssec:stalk}
Let $X$ be a noetherian scheme. For any $x \in X$ the natural map
$\lambda_x \colon \Spec(\OXx) \to X$ is a localizing immersion.
Via the exact global-sections functor~$\Gamma$ on $\Spec(\OXx)$, 
and using~$\Gamma$ to denote~$\R\Gamma$, 
for any $\F \in \Dqct^+(X)$ 
we obtain natural isomorphisms
\[
\Gamma\lambda_x^!\F = \Gamma\lambda_x^*\F 
\cong \F_x.
\]
More generally, if $f \colon X \to Y$ is a map in~$\Se$ and 
$f_x \colon \Spec(\OXx) \to Y$ is the canonical map, then 
for any $\F \in \Dqct^+(Y)$ we obtain natural isomorphisms 
\[
\Gamma f_x^!\F = \Gamma(f\lambda_x)^!\F \cong 
\Gamma\lambda_x^!f^!\F \cong (f^!\F)_x.
\]

Our next application concerns perfect complexes and perfect
maps in~$\Se$. Let us recall some basic facts about them.
We follow mostly the treatment given in Lipman's notes
\cite[4.9]{Lip-LNM}.

\sss
\label{ssec:perfectdef}
Let $X$ be a noetherian scheme. A
complex $\F \in \D(X)$ is called \emph{perfect} if it is locally 
$\D$-isomorphic to a strictly perfect complex, i.e., a 
bounded complex of finite-rank free $\OX$-modules. 
In particular, $\F \in \Dc^{\textup{b}}(X)$.
Let $f \colon X \to Y$ be a scheme-map in~$\Se$.
A complex $\F \in \D(X)$ is called \emph{$f$-perfect}
if it has coherent homology and if it has finite flat $f$-amplitude,
i.e., there exist integers $m \le n$ such that for any $x \in X$,
$\F_x$ is isomorphic in $\D(\cO_{Y,f(x)})$ to a complex
of flat $\cO_{Y,f(x)}$-modules that lives between degrees~$m$ 
and~$n$. If $X = Y$ and $f= 1_X$, then $\F$ is $f$-perfect
iff it is perfect. We call $f$ \emph{essentially perfect}
if~$\OX$ is $f$-perfect.

\sss
\label{ssec:!!to!}
We shall soon give a description
of $f^!$ for an essentially perfect map~$f$ in terms of the relative
dualizing complex~$f^!\OY$, namely a natural isomorphism 
\smash{$f^!(-) \cong f^!\OY \Otimes_{\!X}\, \bL f^*(-)$}. 
Following Lipman's notes, we first define the map 
underlying this isomorphism for any~$f$ in~$\Se$.

Let $X \xto{\;i\;} \ov{X} \xto{\;h\;} Y$ be a factorization of a 
scheme-map $f \colon X \to Y$ in~$\Se$, 
with~$i$ a localizing immersion and $h$ proper. 
Since~$h$ is proper, $h^!$ is right adjoint to~$\R h_*$.
For any 
$\F \in \Dqc^+(Y)$, there results a natural map
\[
\chi^h_{\F} \colon h^!\OY \Otimes_{\!\ov{X}}\, \bL h^*\F 
\to h^! \F,
\]
namely, $\chi^h_{\F}$ is the map adjoint to the composition 
\[
\R h_* (h^!\OY \Otimes_{\!\ov{X}}\, \bL h^*\F) \xto{\; p \;}
\R h_* h^!\OY \Otimes_{\!\ov{Y}}\, \F \xto{\; \tau \;} \F,
\]
where $p$ results from the projection isomorphism 
\cite[Prop.\ 3.9.4]{Lip-LNM}
and~$\tau$ is the trace 
map corresponding to the right-adjointness of~$h^!$ to~$\R h_*$. 
Thus for any factorization $f = hi$ as above 
and for any $\F \in \Dqc^+(Y)$, we obtain a natural map
\[
\chi^{[i,h]}_{\F} \colon f^!\OY \Otimes_{\!X}\, \bL f^*\F \to f^!\F
\]
via the composition of the following natural maps (where
$i^* = \bL i^*$)
\[
f^!\OY \Otimes_{\!X}\, \bL f^*\F 
\cong i^*h^!\OY \Otimes_{\!X}\, i^*\bL h^*\F
\cong i^*(h^!\OY \Otimes_{\!\ov{X}}\, \bL h^*\F) \;{\xto{\chi^h_{\F}}}\; 
i^*h^!\F.
\]

\begin{aprop}
\label{prop:welldefined}
With notation as above, if $f = h_1i_1$ is another factorization 
with~$i_1$ a localizing immersion and~$h_1$ proper, 
then for any $\F \in \Dqct^+(Y)$, it
holds that $\chi^{[i,h]}_{\F} = \chi^{[i_1,h_1]}_{\F}$.
\end{aprop}
Henceforth, for $f, \F$ as above we shall denote $\chi^{[i,h]}_{\F}$ 
by $\chi^f_{\F}$.

\begin{proof}
One follows the same steps as in Lipman's notes 
\cite[4.9.2.2]{Lip-LNM}. For the two factorizations of~$f$ 
as above, first one must find a ``dominating'' factorization.
This can be done using the proof of 
Proposition~\ref{prop:dominating} or one can use this Proposition
to reduce to the finite-type case and then proceed as 
in \emph{loc.\ cit}. The rest of the proof goes through 
using localizing immersions in place of open ones and using 
the properties in \ref{ssec:props}.
\end{proof}

\begin{athm}
\label{thm:perfect}
Let $f \colon X \to Y$ be a map in~$\Se$.
Then the following conditions are equivalent.
\begin{enumerate}
\item The map $f$ is essentially perfect, i.e., $\OX$ is $f$-perfect. 
\item For any open $U \subset X$ and any factorization of
$f|_U$ as \smash{$U \xto{j} Z \xto{g} Y$} where $j$ is a 
closed immersion and $g$ essentially smooth, 
$j_*\OU$ is a perfect $\OZ$-complex. 
\item The complex $f^!\OY$ is $f$-perfect.
\item $f^!\OY \in \Dc^{\textup{b}}(X)$, and for every $\F \in \Dqc^+(Y)$,
the map
\[
\chi^f_{\OY,\F} \colon 
f^!\OY \Otimes_{\!X}\, \bL f^*\F \iso f^!\F
\] 
is an isomorphism.
\end{enumerate}
\end{athm}
\begin{proof}
Here again, one argues as in Lipman's notes 
\cite[Theorem 4.9.4]{Lip-LNM}. Locally on~$X$, the factorization
of~$f$ as in~(ii) always exists and so the proof of loc.\ cit. goes through without any difficulties.
\end{proof}

\sss
\label{ssec:CM-Gor}
Recall that a map $f \colon X \to Y$
is called \emph{Cohen-Macaulay} (resp.~\emph{Gorenstein}) 
if~$f$ is flat and the local rings of the fibers of~$f$ are 
Cohen-Macaulay (resp.~Gorenstein).
Before proceeding further let us recall a well-known fact usually 
not stated in the generality that we need below:
\emph{Let $f \colon X \to Y$ be a flat map in~$\Se$.
Then~$f$ is Cohen-Macaulay iff 
for any connected open subset $U \subset X$, $f^!\OY$ has exactly
one non-vanishing homology and this homology 
is $f$-flat. In particular, $f$ is Gorenstein iff $f^!\OY$ 
is invertible, i.e., on every connected component
of~$X$, the unique nonvanishing homology of~$f^!\OY$ 
is an invertible $\OX$-module.} 

To prove this, first we recall that if $f$ is Cohen-Macaulay,
then on over connected components of~$X$, 
$f$ has a constant relative dimension.
Indeed, for any fiber, each of its connected components
is equidimensional since it is Cohen-Macaulay.
Hence by \cite[15.4.3]{EGAIV}, 
for any integer~$r$, the set of points $x \in X$ 
such that the fiber through~$x$ has dimension~$r$ 
forms an open set. 
Therefore, over any connected open 
subset of~$X$, all the fibers have the same dimension. 
Now we refer to \cite[Theorem 3.5.1]{Con1} or 
\cite[Lemma 1]{Lip1} to complete the proof keeping in 
mind that as the assertions of the previous paragraph 
are local on~$X$ and~$Y$,
we may reduce to the case where~$f$ factors as a 
closed immersion into a formally smooth map. 
Note that the unique non-vanishing 
homology occurs in degree~$-n$ where~$n$ is the (local)
relative dimension. 

In the following proposition we shall use the notation used in 
\ref{ssec:stalk} above, namely that 
for any $x \in X$, $\lambda_x$ denotes the
natural map $\Spec(\OXx) \to X$ while if $f \colon X \to Y$ is a
scheme map, then $f_x$ denotes the natural map 
$\Spec(\OXx) \to Y$.

%
%

\begin{aprop}
\label{prop:properties}
Let $\sfP$ be a property of scheme maps in $\Se$ such that 
for any map $f \colon X \to Y$ in~$\Se$ the following conditions
hold.
\begin{enumerate}
\item[$(\alpha)$]
For any localizing immersion $i \colon W \to X$, if~$f$ 
satisfies~$\sfP$ then so does~$fi$. (We also say that~$\sfP$ 
is stable under localization.)
\item [$(\beta)$]
If $f$ is of \emph{finite type}, then we have the following:
\begin{enumerate}
\item[(i)]
If $f_x$ satisfies $\sfP$ for every $x \in X$, 
then so does~$f$;
\item[(ii)]
the set of all points $x \in X$ such that $f_x$ satisfies $\sfP$ is open.
\end{enumerate}
\end{enumerate}
Then~\textup{(i), (ii)} of~$(\beta)$
hold in general %
and for any~$f$ that satisfies~$\sfP$, there exists a factorization
\smash{$X \xto{\;i\;} Z \xto{\;h\;} Y$} where~$i$ is a localizing 
immersion and~$h$ is a separated finite-type map 
satisfying~$\sf{P}$.
\end{aprop}

\begin{proof}
By Theorem \ref{thm:factorization}, any $\Se$-map 
$f \colon X \to Y$ factors
as $X \xto{\;i\;} Z \xto{\;h\;} Y$ with~$i$ a localizing immersion
and $h$ a separated finite-type map.
For any $x \in X$, with $z = i(x)$ we have $\OXx \cong \OZz$.
Therefore the $\sfP$-locus of~$f$, namely the set of points
$x \in X$ such that $f_x$ satisfies~$\sfP$, 
equals $i^{-1}V$ where~$V$ is the $\sfP$-locus of~$h$,
whence it is open. 

Suppose~$f_x$ satisfies~$\sfP$ for each $x \in X$. Then~$V$ 
contains~$i(X)$ and therefore, using~$V$ in place of~$Z$ 
and replacing $i,h$ suitably we may assume that in the 
factorization $f= hi$ above, $h$ satisfies~$\sfP$. Thus~$f$ 
satisfies~$\sfP$. Conversely if~$f$ satisfies~$\sfP$, then
so does~$f_x$ and hence arguing as above yields the desired
factorization of~$f$. 
\end{proof}

\sss
\label{ssec:P-examples}
\textsc{Examples}
\newline
We give some examples of properties where 
Proposition~\ref{prop:properties} applies. A minor note here
is that as a local condition, to have $f_x$ satisfy $\sfP$
could a-priori be viewed as a strong requirement 
since the presence 
of the non-closed points in~$\Spec(\OXx)$ means 
the behavior of~$f$ at those points too is being considered. 
However in the presence of the finiteness hypothesis on~$f$, 
it remains equivalent to more familiar local versions of~$\sfP$.

In what follows, each of the properties considered is stable
under localization, so we only discuss how (i) and~(ii) 
of Proposition~\ref{prop:properties}($\beta$) hold.

\smallskip

(a). $\sfP=$ Flat. Here $f_x$ being flat is equivalent to~$f$
being flat at~$x$ in the usual sense, i.e., $\OXx$ is a flat
$\OYy$-module. Thus  (i) holds and for (ii) we refer 
to \cite[11.1.1]{EGAIV}.

(b). $\sfP=$ Essentially Perfect. As in (a), $f_x$ being essentially  
perfect is equivalent to $\OXx$ being a perfect $\OYy$-complex.
Thus~(i) holds while for~(ii) we note that in the notation 
of Theorem~\ref{thm:perfect}, 
the perfect locus of~$f|_U$ equals the 
inverse image under~$j$ of the perfect locus of the coherent 
$\OZ$-module $j_*\OU$, which is open.

(c). $\sfP=$ Cohen-Macaulay (CM). Since the fibers of $f_x$ 
are again localizations of the fibers of~$f$ and~(i) holds for 
flatness we see that~(i) holds in this case too. For (ii) 
we argue using \cite[12.1.6]{EGAIV}. Indeed, by \emph{loc. cit}
the set~$U$ of points $x \in X$ such that the fiber through~$x$
has a CM local ring at~$x$ forms an open set.
If~$f_x$ is CM, then $x \in U$ and conversely
if~$x$ is in~$U$, then the natural image of $\Spec(\OXx)$ in~$X$
is contained in~$U$, whence $f_x$ is CM. 
Thus (ii) holds.

(d). $\sfP=$ Gorenstein. For (i) we argue as in (c). For~(ii), 
first, by~(a) above, we may shrink~$X$ if necessary and
assume that~$f$ is flat.
If~$f_x$ is Gorenstein, then using the
isomorphism $(f^!\OY)_x \cong \Gamma f_x^!\OY$
(see \ref{ssec:stalk}) and exactness of~$\Gamma$ we obtain
from \ref{ssec:CM-Gor} that $(f^!\OY)_x$ has a unique 
homology which is a free $\OXx$-module of rank 1. Since 
$f^!\OY \in \Dc^{\textup{b}}(X)$, it follows that~$x$ has an open
neighborhood~$U$ such that $(f^!\OY)|_U$
has a unique invertible homology, whence by \ref{ssec:CM-Gor}, 
(ii) holds.

(e). $\sfP=$ Essentially smooth. If~$f_x$ is (essentially) smooth 
then it is also flat. Thus smoothness for each $f_x$ implies
flatness for~$f$ and moreover that every fiber of~$f$, being locally smooth, is smooth. Since flatness + smooth fibers is equivalent
to smoothness (\cite[17.5.1]{EGAIV}), (i) holds. 
For~(ii), we argue as in (c).

(f). $\sfP=$ Unramified. By \cite[Prop 2.6.1]{LNS}, 
the module of relative differentials 
$\Omega_{f_x}$ on $\Spec(\OXx)$ is coherent. 
Moreover, there is a natural isomorphism 
$\Omega_{f_x} \cong \lambda_x^*\Omega_{f}$ 
as the natural map 
$\lambda_x \colon \Spec(\OXx) \to X$ is (essentially) \'{e}tale. 
As in \ref{ssec:stalk}, there result natural isomorphisms
$\Gamma\Omega_{f_x} \cong \Gamma\lambda_x^*\Omega_{f}
\cong (\Omega_f)_x$ where~$\Gamma$ is the 
(exact and faithful) 
global-sections functor on $\Spec(\OXx)$. Thus~(i) and~(ii) 
can be verified by looking at stalks of~$\Omega_f$.
 
(g). $\sfP=$ Essentially \'{E}tale. 
One can use (e) while keeping track of 
the rank of module of differentials as in (f). 

\medskip

\textsc{Acknowledgements}.
I thank Professor Joseph Lipman for suggesting the main
compactification problem to me and for various discussions on it.

\providecommand{\bysame}{\leavevmode\hbox to3em{\hrulefill}\thinspace}
\providecommand{\MR}{\relax\ifhmode\unskip\space\fi MR }
\providecommand{\MRhref}[2]{%
  \href{http://www.ams.org/mathscinet-getitem?mr=#1}{#2}
}
\providecommand{\href}[2]{#2}

\end{document}